# Fuzzy set representation of a prior distribution

## Glen Meeden[*,1]

*University of Minnesota*

**Abstract:** In the subjective Bayesian approach uncertainty is described by a prior distribution chosen by the statistician. Fuzzy set theory is another way of representing uncertainty. Here we give a decision theoretic approach which allows a Bayesian to convert their prior distribution into a fuzzy set membership function. This yields a formal relationship between these two different methods of expressing uncertainty.

## Contents

## 1. Introduction

For a subjective Bayesian uncertainty about the unknown parameter or state of nature can be expressed through a prior distribution. If $\theta$ denotes a typical parameter value and $\Theta$ the set of all possible parameter values then the prior distribution over $\Theta$ summarizes their knowledge and beliefs about the parameter.

Fuzzy set theory, introduced in Zadeh [7], is another approach to representing uncertainty. A fuzzy set $A$, a subset of $\Theta$, is characterized by its membership function. This is a function defined on $\Theta$ whose range is contained in the unit interval. At a point $\theta$ the value of the membership function is a measure of how much we think $\theta$ belongs to the set $A$. Statisticians have been slow to embrace fuzzy set theory. Taheri [6] gives a review of applications of fuzzy set theory concepts to statistical methodology. Bayesians have shown less interest in fuzzy ideas than frequentists. Singpurwalla and Booker [5] have proposed a model which incorporates fuzzy membership functions into a subjective Bayesian setup. However, they do not give membership functions a probabilistic interpretation. In the imprecise or vague approach to Bayesian statistics a decision maker selects a family of possible prior

*Supported in part by NSF Grant DMS-04-06169.
[1]School of Statistics, University of Minnesota, 313 Ford Hall, 224 Church Street, SE, Minneapolis, MN 55455, USA, e-mail: glen@stat.umn.edu

distributions to represent their prior beliefs. de Cooman [2] presents an uncertainty model for vague probability assessments that is closely related to Zadeh's approach [7].

The concept of confidence intervals is a frequentist approach to expressing uncertainty about an unknown parameter given data. It has long been recognized that naive users have difficulty interpreting confidence intervals. They have a tendency to give a probabilistic interpretation to the observed confidence interval.

It has also long been known that for discrete data conventional confidence intervals, which we will also call "crisp" confidence intervals, using a term from fuzzy set theory, can perform poorly. A recent article [1] reviews the problems with crisp confidence intervals for binomial models. Because of the inherent flaws in crisp confidence intervals for discrete problems a new confidence interval notion has been suggested called fuzzy confidence intervals [3]. Given the data a fuzzy confidence interval is just the membership function of the set of plausible or reasonable values for $\theta$. One way to think about such membership functions is that they are generalizations of randomized intervals where no randomization is ever implemented. They argued that fuzzy confidence intervals overcome the difficulties of the usual crisp intervals for discrete probability models.

In terms of frequency of coverage discrete data Bayesian credible intervals will suffer from the same problem that conventional intervals do. This should be of concern to objective Bayesians who want their intervals to have good frequentist properties. One way to approach this problem is to find a method that allows them to use their posterior to get a sensible fuzzy interval instead of the usual Bayesian credible interval.

Here we consider a no data statistical decision problem where the set of possible decisions is the class of all membership functions defined on $\Theta$. We then define a family of loss functions. These functions measure the loss incurred when a probability distribution is replaced by a fuzzy membership function. For any loss function in the family and a given prior distribution we solve the resulting no data decision problem. This gives a method for converting a prior or posterior into a fuzzy membership function. For a given fuzzy membership function we also study the problem of identifying the family of prior distributions whose common solution to the no data decision problem is this function. This sets up a formal relationship between the two theories.

## 2. Fuzzy set theory

We will only use some of the basic concepts and terminology of fuzzy set theory, which can be found in the most elementary of introductions to the subject [4].

A *fuzzy set* $A$ in a space $\Theta$ is characterized by its *membership function*, which is a map $I_A : \Theta \to [0,1]$. The value $I_A(\theta)$ is the "degree of membership" of the point $\theta$ in the fuzzy set $A$ or the "degree of compatibility ... with the concept represented by the fuzzy set." See ([4], p. 75). The idea is that we are uncertain about whether $\theta$ is in or out of the set $A$. The value $I_A(\theta)$ represents how much we think $\theta$ is in the fuzzy set $A$. The closer $I_A(\theta)$ is to 1.0, the more we think $\theta$ is in $A$. The closer $I_A(\theta)$ is to 0.0, the more we think $\theta$ is not in $A$.

A fuzzy set whose membership function only takes on the values zero or one is called *crisp*. For a crisp set, the membership function $I_A$ is the same thing as the indicator function of an ordinary set $A$. Thus "crisp" is just the fuzzy set theory way of saying "ordinary," and "membership function" is the fuzzy set theory way of saying "indicator function." The *complement* of a fuzzy set $A$ having membership function $I_A$ is the fuzzy set $B$ having membership function $I_B = 1 - I_A$.

If $I_A$ is the membership function of a fuzzy set $A$, the $\gamma$-*cut* of $A$ ([4], Section 5.1) is the crisp set
$$^\gamma I_A = \{\theta : I_A(\theta) \geq \gamma\}.$$

Clearly, knowing all the $\gamma$-cuts for $0 \leq \gamma \leq 1$ tells us everything there is to know about the fuzzy set $A$. The 1-cut is also called the *core* of $A$, denoted core($A$) and the set
$$\text{supp}(A) = \bigcup_{\gamma > 0} {}^\gamma I_A = \{\theta : I_A(\theta) > 0\}$$
is called the *support* of $A$ ([4], p. 100). A fuzzy set is said to be *convex* if each $\gamma$-cut is convex ([4], pp. 104–105).

## 3. A decision problem

For simplicity we assume that $\Theta$ is an interval of real numbers and the prior $\pi$ is a continuous probability density function defined on it.

Let $\mathcal{A}$ be the class of all measurable membership functions defined on $\Theta$. Then $\mathcal{A}$ is the space of possible decisions or actions with a typical member denoted by $A$. Given a prior density $\pi$ on $\Theta$ we want to find the membership function or fuzzy set $A$ which best represents $\pi$. We do this by defining a loss function and then solving the no data decision problem.

Our loss function will depend on four known parameters which are specified by the statistician. They are $a_1 \geq 0$, $a_2 \geq 0$, $b_1 \geq 0$ and $b_2 \geq 0$ where at least one of the $a_i$'s and at least one of the $b_i$'s must be strictly positive. Then the loss incurred when action $A$ is taken and $\theta$ is the true state of nature is given by

$$(1) \quad L(A, \theta) = a_1\{1 - I_A(\theta)\} + \frac{a_2}{2}\{1 - I_A(\theta)\}^2 + \int_\Theta \left\{b_1 I_A(\theta) + \frac{b_2}{2}(I_A(\theta))^2\right\} d\theta.$$

To understand this loss function remember that we want to find the fuzzy set or membership function $A$ which best represents the set of sensible or reasonable parameter values under our prior $\pi$. Hence if $\theta$ is the true parameter point we want $I_A(\theta)$ to be close to 1. This explains the presence of the first two terms in equation (1). But on the other hand we do not want the fuzzy set to be too large. This is controlled by the last term in the equation which is a measure of the overall size of the fuzzy set.

We now find the solution for this no data decision problem.

**Theorem 1.** *Let $\pi(\theta)$ be a prior density on $\Theta$. Then for the loss function of equation (1) the fuzzy set membership $A$ which satisfies*

$$\int_\Theta L(A, \theta) \pi(\theta) \, d\theta = \inf_{A' \in \mathcal{A}} \int_\Theta L(A', \theta) \pi(\theta) \, d\theta$$

*is given by*

$$(2) \quad I_A(\theta) = \begin{cases} 0, & \text{for} \quad 0 \leq \pi(\theta) < b_1/(a_1 + a_2), \\ \frac{(a_1 + a_2)\pi(\theta) - b_1}{a_2 \pi(\theta) + b_2}, & \text{for} \quad b_1/(a_1 + a_2) \leq \pi(\theta) \leq (b_1 + b_2)/a_1, \\ 1, & \text{for} \quad \pi(\theta) > (b_1 + b_2)/a_1. \end{cases}$$

*Proof.* Note that we can write

$$\int_\Theta L(A',\theta)\pi(\theta)\,d\theta = \int_\Theta \left\{ a_1\{1 - I_{A'}(\theta)\} + \frac{a_2}{2}\{1 - I_{A'}(\theta)\}^2 \right\}\pi(\theta) + \\ b_1 I_{A'}(\theta) + \frac{b_2}{2}(I_{A'}(\theta))^2 \right\} d\theta$$

so that to find the solution it is enough to minimize the integrand of the previous equation for each fixed value of $\theta$. But for a fixed $\theta$ the integrand is just a quadratic function of $I_{A'}(\theta)$ and a simple calculus argument completes the proof. □

The theorem remains true when $a_1 = 0$ if we assume dividing by zero yields infinity.

Note that the solution is unchanged if the loss function is multiplied by a positive number. Without loss of generality we could set one of the four parameters defining the loss function equal to one but having four parameters will be convenient in the following discussion.

As with any decision problem the solution depends strongly on the loss function. We believe our family of loss functions is flexible and captures some of the important aspects of the problem. Finding a good fuzzy set to summarize our information about a parameter is much like finding a good credible set. We want it to include the likely values but without it getting to large. The loss function in equation (1) is essentially the sum of two quadratic functions. The first part is quadratic in non-membership in the set of likely values while the second part is quadratic in a measure of the size of the set. If we just include the linear terms in each part then the optimal solution will always be a crisp set. It is necessary to include the quadratic terms to get a true fuzzy set as a solution.

We see from equation (2) that the optimal membership function is related to the prior $\pi$ in a sensible fashion. The solution is 1 where the prior is large, 0 where the prior is small and a rescaling between the two cases. Note that for a given bounded $\pi$ if $b_1$ is chosen large enough then the solution to our decision problem is the membership function which is identically zero. On the other hand if $\pi$ is bounded away from zero and $a_1$ is chosen large enough then the solution to our decision problem is the membership function which is identically one.

## 4. Relating priors and fuzzy sets

We have considered the problem of converting a prior distribution into a fuzzy membership function. In some situations it could be of interest to be able to move in the other direction. That is, transform the uncertainty expressed in a fuzzy membership function into the Bayesian paradigm. One way to do this would be to find a loss function and prior for which the solution to our decision problem is the fuzzy membership function in hand. This suggests the following three questions.

- For a specified fuzzy membership function, $I_A$, and a specified loss function does there exist a prior density function for which the solution to our decision problem is $I_A$?
- For a specified fuzzy membership function, $I_A$, does there exist a loss function and a prior density function for which the solution to our decision problem is $I_A$?
- If a solution does exist for question 1 is it unique?

We see from equation (2) that for $I_A$ to be a solution for $\pi$ we must have

$$\pi(\theta) = \frac{b_1 + b_2 I_A(\theta)}{a_1 + a_2(1 - I_A(\theta))} \quad \text{for } \theta \text{ where} \quad 0 < I_A(\theta) < 1 \tag{3}$$

From this we see that the answer to our first question is no. This is because when $\Theta$ is unbounded $\pi$ in the previous equation need not be integrable and even when it is it need not integrate to one. The answer to the second question is yes whenever $I_A(\theta)$ is integrable. Since in this case we can always select $b_1 \geq 0$ and $b_2 > 0$ to make $\pi(\theta)$ of equation (3) a density. When a solution exists it need not be unique.

For a simple example we set $a_2 = 0$ and let the other three parameters be positive. Consider the special case where $\Theta$ is bounded. If we set

$$r_1 = b_1/a_1 \quad \text{and} \quad r_2 = (b_1 + b_2)/a_1 \tag{4}$$

we find that

$$a_1 = b_2/(r_2 - r_1) \quad \text{and} \quad b_1 = r_1/(r_2 - r_1) \tag{5}$$

and the solution from equation (2) has the form

$$I_A(\theta) = \begin{cases} 0, & \text{for } 0 \leq \pi(\theta) < r_1, \\ (\pi(\theta) - r_1)/(r_2 - r_1), & \text{for } r_1 \leq \pi(\theta) \leq r_2, \\ 1, & \text{for } \pi(\theta) > r_2. \end{cases} \tag{6}$$

Now let $I_A$ be given and assume that the length of $\Theta$ is $\ell$. If $r_1 < 1/\ell$ then there exist a unique $r_2 > r_1$ such that

$$\pi_{A,r_1}(\theta) = (r_2 - r_1)I_A(\theta) + r_1 \quad \text{for } \theta \in \Theta \tag{7}$$

is a prior distribution over $\Theta$. Moreover we can find values for $a_1$, $b_1$ and $b_2$ which satisfy equation (4). With this loss function $I_A$ will be the solution to our decision problem when the prior is $\pi_{A,r_1}$. Furthermore if the sets where $I_A(\theta) = 0$ and $I_A(\theta) = 1$ each have positive Lebesgue measure then it will not be the unique prior with this property. Any prior density $\pi$ satisfying

$$\begin{aligned} \pi(\theta) &\leq r_1 \quad \text{when} \quad I_A(\theta) = 0, \\ \pi(\theta) &= \pi_{A,r_1}(\theta) \quad \text{when} \quad 0 < I_A(\theta) < 1, \\ \pi(\theta) &\geq r_2 \quad \text{when} \quad I_A(\theta) = 1 \end{aligned} \tag{8}$$

will also be a solution for our decision problem.

Among the set of possible solutions the one in equation (7) has two nice properties. First of all it is continuous whenever $I_A(\theta)$ is continuous. Secondly it treats the members of $\{\theta : I_A(\theta) = 1\}$ similarly and the members of $\{\theta : I_A(\theta) = 0\}$ similarly. More importantly, this identification of a fuzzy membership function with a class of prior distributions demonstrates that we can give roughly equivalent expressions of uncertainty in the Bayesian and fuzzy paradigms.

Finally, we address the question of uniqueness. The previous discussion indicates that if we want uniqueness we should consider membership functions which never take on zero or one as a possible value. Let $I_A$ be such a membership function and let $a_1 > 0$ and $a_2 > 0$ be fixed and suppose $\Theta$ is the unit interval. Then integrating

equation (3) we have

$$\int_0^1 \pi(\theta)\,d\theta = b_1 \int_0^1 \frac{1}{a_1 + a_2(1 - I_A(\theta))}\,d\theta + b_2 \int_0^1 \frac{I_A(\theta)}{a_1 + a_2(1 - I_A(\theta))}\,d\theta$$
$$= b_1 c_1(a_1, a_2) + b_2 c_2(a_1, a_2).$$

Hence $\pi$ will be a probability density function whenever

$$b_1 \in [0, 1/c_1(a_1, a_2)] \quad \text{and} \quad b_2 = (1 - b_1 c_1(a_1, a_2))/c_2(a_1, a_2).$$

To better understand the relationship between $I_A$ and its corresponding prior we consider a simple example. Let

(9) $$I_A(\theta) = 6.075\,\theta^2(1 - \theta) \quad \text{for } \theta \in [0, 1].$$

We consider two different choices of the $a_i$'s and for each case two different choices of $b_1$. For the first case $a_1 = 1$ and $a_2 = 7$. The maximum possible value for $b_1$ is 3.40 and our two choices for the $b_i$'s are $b_1 = 0.01, b_2 = 5.15$ and $b_1 = 3.35, b_2 = 0.072$. In the second case $a_1 = 4$ and $a_2 = 2$. The maximum possible value for $b_1$ is 4.91 and our two choices for the $b_i$'s are $b_1 = 0.01, b_2 = 9.02$ and $b_1 = 4.50, b_2 = 0.76$. For each of the four combinations we found the unique prior whose solution to the decision problem yields the fuzzy membership function of equation (9). The membership function along with the four priors are shown in the figure.

The membership function is the solid curve. The two curves with the two largest maximums are the solutions for the first case where $a_1 = 1$ and $a_2 = 7$. Of the two solutions the one with $b_1 = 0.01$ has the largest maximum. The other two curves are the solutions for the second case. Again the solution for $b_1 = 0.01$ has the largest of the two maximums. These curves demonstrate what a closer inspection of equation (3) yields. For a fixed $a_1$ and $a_2$ the solution becomes less concentrated about its mode as $b_1$ increases from zero to its maximum value. Also the solution becomes less concentrated about its mode as we increase $a_1$ and decrease $a_2$. But in all cases the priors do reflect the shape of their common solution.

An interesting consequence of this unique correspondence is that it gives a way to update a large class of fuzzy membership functions given data. Suppose an expert has selected a fuzzy membership function to represent their uncertainty. The statistician then selects appropriate values of the $a_i$'s and the $b_i$'s and uses equation (3) to transform it into a prior. Then given the data they find the posterior distribution which is then converted back to a fuzzy membership function using the theorem with the $a_i$ and $b_i$ values.

This result is somewhat surprising since a fuzzy membership function must satisfy less conditions then a probability density function since it need not be integrable. At first glance the previous example where a membership function corresponded to a family of priors seems more reasonable. To get the unique correspondence, however, we made two fairly strict assumptions. The function in equation (3) needed to be integrable and the range of the membership function had to lay in the open unit interval. Both these conditions on the membership function seem not so surprising if we hope to convert it to a probability density function.

## 5. Some final remarks

Mainline statistics has shown little interest in fuzzy set theory. This is especially true for most Bayesians since they believe that they already have a good way to express uncertainty. Here we have argued that Bayesians should be more interested in fuzzy set theory. For discrete data, just as for frequentists, there are certain advantages to considering interval estimates as fuzzy sets. We noted that our scheme for converting a prior density into a fuzzy membership function could also be used to relate some fuzzy membership functions to prior densities. In some cases a fuzzy membership function will correspond to a family of densities while under more restricted conditions it will correspond to a unique density. The relationship seems intuitively sensible and as far as we know it is the first simple formal correspondence between the two theories which until now have lived in different worlds.

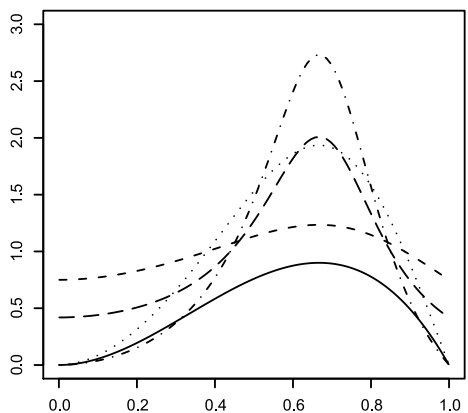

Fig 1. *A plot of the fuzzy membership function (the solid line) in equation* (9) *and four priors whose common solution for four loss functions is the fuzzy membership function. The two priors for the $a_1 = 1$ and $a_2 = 7$ case are the ones with the two largest maximums. The other two priors are for the two $a_1 = 4$ and $a_2 = 2$ cases.*

A copy of Geyer and Meeden [3] and related material can be found at http//:www.stat.umn.edu/~glen/papers/.

**Acknowledgments.** The author would like to thank a referee for a helpful comment.